\documentclass{amsproc}
\usepackage{cite}
\newtheorem{theorem}{Theorem}[section]
\newtheorem{lemma}[theorem]{Lemma}
\newtheorem{prob}[theorem]{Problem}

\newtheorem{cor}[theorem]{Corollary}
\newtheorem{prop}[theorem]{Proposition}
\theoremstyle{definition}

\newtheorem{example}[theorem]{Example}

\theoremstyle{remark}

\numberwithin{equation}{section}

\begin{document}

\title[DERIVATIONS, LOCAL AND 2-LOCAL DERIVATIONS ON ALGEBRAS]{DERIVATIONS,
LOCAL AND 2-LOCAL DERIVATIONS ON ALGEBRAS OF MEASURABLE OPERATORS}

\author{Shavkat Ayupov}
\address{Dormon yoli 29,
Institute of
 Mathematics,
 National University of
Uzbekistan,
 100125  Tashkent,   Uzbekistan}
\email{sh$_{-}$ayupov@mail.ru}

\thanks{The authors would like to acknowledge the hospitality
of the California State University, Fullerton, during
USA-Uzbekistan Conference on Analysis and Mathematical Physics,
May 20-23, 2014.}

\author{Karimbergen Kudaybergenov}
\address{Ch. Abdirov 1,
Department of Mathematics,
 Karakalpak state university,
 Nukus
230113, Uzbekistan}
\email{karim2006@mail.ru}


\date{January 1, 1994 and, in revised form, June 22, 1994.}



\begin{abstract}
The present paper presents a survey of some recent results devoted to
 derivations, local derivations and 2-local derivations on various algebras of measurable operators
 affiliated with von Neumann algebras. We give a complete description of derivation on these algebras,
 except  the case where the von Neumann algebra is of type II$_1$. In the latter case the result is obtained under an extra
  condition of measure continuity of derivations. Local and 2-local derivations
 on the above algebras are also considered. We give sufficient conditions on a von Neumann algebra $M$,
 under which every local or 2-local derivation on the algebra of measurable operators  affiliated with $M$
 is automatically becomes a derivation. We also give examples of commutative algebras of measurable
 operators admitting local and 2-local derivations which are not derivations.
\end{abstract}

\maketitle

\tableofcontents

\section{Introduction}

The present paper is devoted to some recent results
concerning derivations and derivation-type mappings on certain classes of unbounded operator
algebras.

The theory of algebras of operators acting on a Hilbert space began
in 1930s with a series of papers by Murray and von Neumann  (see
\cite{MN1, MN2, MN3, MN4}), motivated by  the theory of unitary
group representations and certain aspects of the quantum mechanical formalism.
They analyzed the structure of the family of algebras which are
referred nowadays as von Neumann algebras or $W^\ast$-algebras and which
have the distinctive property of being closed in
the weak operator topology. In 1943 Gelfand and Naimark
developed the theory of uniformly closed operator $^\ast$-algebras,
which are now called $C^\ast$-algebras.

Nowadays the theory of operator algebras  plays an important role both in pure
mathematical and application aspects. This is motivated by the fact that in terms of  operator
algebras, their states, representations, groups of automorphisms,
and derivations one can describe and investigate  properties of
model systems  in the quantum field theory and statistical physics.

 Let
$\mathcal{A}$ be an algebra over the field of complex numbers. A
linear (respectively, additive) operator $D:\mathcal{A}\rightarrow \mathcal{A}$
is called a linear (respectively, additive) \textit{derivation} if it satisfies
the identity $D(xy)=D(x)y+xD(y)$ for all  $x, y\in \mathcal{A}$
(Leibniz rule). Each element  $a\in \mathcal{A}$ defines a linear
derivation $D_a$ on $\mathcal{A}$ given by $D_a(x)=ax-xa,\,x\in
\mathcal{A}.$ Such derivations $D_a$ are said to be \textit{inner}.
If the element  $a$ implementing the derivation
$D_a$ on $\mathcal{A},$ belongs to a larger algebra $\mathcal{B},$
containing  $\mathcal{A}$ (as a proper ideal as usual) then $D_a$
is called a  \textit{spatial derivation}.

One of the main problems considered in the theory of derivations
is to prove the automatic continuity, innerness or spatialness of
derivations, or to show the existence of non inner and
discontinuous derivations on various topological algebras.

In particular, it is a general algebraic problem to find  algebras
which admit only inner derivations.

A more general problem is the following one:  given an algebra
$\mathcal{A},$ does there exist an algebra $\mathcal{B}$ containing $\mathcal{A},$
such that any derivation of the algebra $\mathcal{A}$ is spatial and
 implemented by an element from $\mathcal{B} ?$ (see e.g.~\cite{Che}, \cite{Sak1}).

 The theory of derivations in operator algebras is an important and well investigated part of the general theory
 of operator algebras, with applications in mathematical physics (see, e.g.~\cite{Bra},
\cite{Sak1}, \cite{Sak2}). It is well known that every derivation
of a
  $C^{\ast}$-algebra is bounded (i.e. is norm continuous), and that every derivation of a von
 Neumann algebra is  inner.
 For  a detailed exposition of the theory of bounded derivations we refer to the monographs of
 Sakai \cite{Sak1}, \cite{Sak2}. A comprehensive study of derivations in general
 Banach algebras is given in the monograph of Dales~\cite{Dal} devoted to the study of
 automatic continuity of derivations on various classes of Banach
 algebras.

 Investigations of general unbounded derivations (and derivations on unbounded operator algebras)
 began much later and were motivated mainly by needs of mathematical physics, in particular by the
 problem of constructing the dynamics in quantum statistical mechanics. The kinematical structure
 of a physical system in the quantum field theory (systems with  infinite number of degrees
of freedom) is described by an operator algebra $A$, where states are positive normalized linear functionals
on  $A$, and observables are elements of this algebra $A$. The dynamical evolution of the system
is given by a group of $\ast$-automorphism of the operator algebra $A.$
The infinitesimal motion is described by some form of Hamiltonian
formalism, incorporating the interparticle interaction. In quantum field theory the infinitesimal motion is given
 by a derivation $d$ on the operator algebra $A$ of observables. The basic problem which occurs in this
 approach – is the integration of these infinitesimal motion in order to obtain the
dynamical flow. In terms of operator algebras this means:
to prove that a given derivation on the algebra of observables
is the infinitesimal generator of a one-parameter automorphisms
group, moreover it is spatial (i.e. defined by some Hamiltonian
operator) or even inner (i.e. the Hamiltonian operator is itself
an observable in the considered physical system). For details we refer to \cite{Bra}.

 The development of a non commutative integration theory was
  initiated by Segal~\cite{Seg},
 who considered new classes of (not necessarily Banach) algebras of unbounded operators, in particular the
 algebra  $S(M)$ of all measurable operators affiliated with a von
 Neumann algebra $M.$ Since  the algebraic, order and topological properties of
 the algebra $S(M)$ are somewhat similar to those of $M$,  in~\cite{Ayu1}, \cite{Ayu2005}
 the above problems have been considered for derivations on the algebra $S(M).$
  If the von Neumann algebra $M$ is abelian then it is *-isomorphic to the algebra
   $L^{\infty}(\Omega)=L^{\infty}(\Omega, \Sigma,
\mu)$ of all (classes of equivalence of)  essentially bounded measurable complex functions on a measure space
$(\Omega, \Sigma, \mu)$ and therefore,  $S(M)\cong
L^{0}(\Omega),$   where $L^{0}(\Omega)=L^{0}(\Omega, \Sigma,
\mu)$ is the algebra of all   measurable complex functions on
$(\Omega, \Sigma, \mu),$ and hence, in this case   inner derivations on $S(M)$ are identically zero, i.e. trivial.

In the abelian case Ber, Sukochev, and
Chilin in~\cite{Ber} obtained necessary and sufficient
conditions for existence of non trivial derivations on
commutative regular algebras. In particular they  prove that
the algebra  $L^{0}(0, 1)$ of all  measurable complex functions on
the interval $(0, 1)$ admits non trivial derivations.
Independently, Kusraev (see~\cite{Gut}, \cite{Kus2},
\cite{Kus}) by means of Boolean-valued analysis establishes
necessary and sufficient conditions for existence of non
trivial derivations and automorphisms on extended complete complex
$f$-algebras. In particular, he also proves that the algebra $L^{0}(0, 1)$ admits
non trivial derivations and automorphisms.  It is
clear that these derivations are discontinuous in the measure
topology, and  they are neither inner nor spatial.
Therefore, the properties of derivations on the
algebra $S(M)$ of unbounded operators are very far from being similar to those exhibited
by derivations on
$C^{\ast}$- or von Neumann algebras. But it seems that the
existence of such "exotic" examples of derivations is deeply
connected with the commutativity of the underlying von Neumann
algebra $M.$ In view of this conjecture the present authors
suggested to investigate the above problems in a non commutative
setting (see~\cite{Alb3}, \cite{Alb1}), by considering derivations on
the algebra $LS(M)$ of all locally measurable operators with
respect to a semi-finite von Neumann algebra $M$ and on various
subalgebras of $LS(M).$ The most complete results concerning
derivations on $LS(M)$ have been obtained by the authors and
collaborators in the case of type I von Neumann algebras. Some of our results have been confirmed
 independently in~\cite{Ber1} by representation of measurable operators
 as operator valued functions. Another  approach to similar problems
 in the framework of type I $AW^{*}$-algebras has been outlined in~\cite{Gut}.

 The  paper is organized as follows. In section~\ref{prem} we present the preliminaries and basic results on
 non commutative integration theory and recall  definitions of the algebras  $S(M)$ of
  measurable operators,  $LS(M)$ of locally measurable operators affiliated with
 a von Neumann algebra $M$. We also consider their subalgebras: $S(M,\tau)$ of  $\tau$-measurable
 operators, and $S_0(M, \tau)$ of $\tau$-compact operators affiliated with
 the von Neumann algebra $M$ and  a faithful normal semi-finite trace  $\tau$ on $M.$
 The latter algebras  equipped with the measure topology become metrizable topological algebras.
 Section~\ref{derI} contains a complete description of derivations on the algebras  $LS(M),$
 $S(M)$, $S(M, \tau)$ and $S_0(M, \tau)$ for a type I von Neumann algebra $M.$ We give a general construction
 of  derivations which are neither inner nor spatial, and moreover, which are
 discontinuous in the measure topology on the algebra $LS(M)$ = $S(M)$ for a finite type I von Neumann
 algebra $M$.  We show that for properly infinite type I von Neumann algebra $M$, the algebras
 $LS(M),$ $S(M)$ and $S(M, \tau)$ admit only inner derivations. Derivations on  the algebra $S_0(M, \tau)$ of
 $\tau$-compact operators are investigated for arbitrary semi-finite
 (i.e. type II algebras are also included) von Neumann algebras.
 We show that in the properly infinite case every derivation on this algebra is spatial and implemented
 by an element of $S(M, \tau)$.
 In Section~\ref{derG} we extend the results of the previous section to additive derivations
 on $LS(M)$ for type  I$_\infty$ or type III von Neumann algebras. Here we  also  present some recent results of \cite{Ber2013}
 which generalize this theorem for arbitrary properly
 infinite von Neumann algebras. The problem of description of derivationson $S(M)$  remains open  only when $M$ is of type II$_1.$
 We present a positive solution of this problem  in the case of derivations which are continuous in the measure topology.

 In Section~\ref{derloc} we study the
 so-called local derivations on the algebra $S(M, \tau)$. This notion was introduced
 by Kadison, who investigated such mappings on von Neumann algebras and some
 polynomial algebras. Here we extend his results and show that every continuous
 (in  the measure topology) local derivation on $S(M, \tau)$ is a derivation.
 In the case of an abelian von Neumann algebra $M$ we  give  necessary and
 sufficient conditions for the existence of local derivations on $S(M, \tau)$
 which are not derivations. Finally, in Section~\ref{der2loc} we consider   $2$-local derivations
 on algebras of measurable operators. Such mappings were introduced by Semrl, who obtained their description in the
 case of the algebra $B(H)$ for infinite dimensional separable Hilbert space $H$. Here we give
 the exposition of  results which describe $2$-local derivations on the algebra
 $S(M)$ of measurable operators affiliated with an arbitrary von Neumann algebra $M$ of type I.

\section{Locally measurable operators affiliated with von Neumann algebras}
\label{prem}

Let $H$ be a Hilbert space over the field $\mathbb{C}$ of complex
numbers, and let $B(H)$ be the algebra of all bounded linear
operators on $H.$ Denote by $\mathbf{1}$ the identity operator on
$H$, and let $P(H)=\{p\in B(H): p=p^2=p^{\ast}\}$ be the lattice of
projections in $B(H).$ Consider a von Neumann algebra $M$ on $H,$
i.e. a  *-subalgebra of $B(H)$  closed in the week operator topology and containing the
operator $\mathbf{1}$. Denote by $\|\cdot\|_M$ the operator norm
on $M.$ The set $P(M)=P(H)\cap M$ is a complete orthomodular lattice
with respect to the natural partial order on $M_h=\{x\in M:
x=x^{\ast}\},$ generated by the cone $M_{+}$ of positive operators
from $M.$

Two projections $e,f\in P(M)$ are said to be \textit{equivalent}
(denoted by $e\sim f$) if there exists a partial isometry $v\in M$
with  initial projection $e$ and  final projection $f$, i.e.
$v^{\ast}v=e,\, vv^{\ast}=f $.
The relation $"\sim"$ is  equivalence relation on the lattice
$P(M).$

A projection $e\in P(M)$ is said to be \textit{finite}, if for $f\in
P(M),\, f\leq e,\, f\sim e$ implies that $e=f.$

A von Neumann algebra $M$ is said to be

-- \textit{finite} if $\mathbf{1}$ is a finite projection;

-- \textit{semi-finite} if every non zero projection in $M$
admits a nonzero finite sub-projection;

-- \textit{infinite} if $\mathbf{1}$ is not finite;

-- \textit{properly infinite}, if every non zero central projection
in $M$ is infinite (i.e. not finite);

--\textit{purely infinite or type III} if every non zero projection
in $M$ is infinite.

A projection $e$ in a von Neumann algebra $M$ is said to be abelian
if $eMe$ is an abelian  von Neumann algebra. Since the lattice of projection $P(M)$
is complete, for every projection $e$ in $M$ there exists the least central projection  $z(e)$ containing $e$ as a sub-projection;
it is called the central support of $e$. A projection $e$ is said to be faithful if $z(e) =\mathbf{1}$.
A von Neumann algebra $M$ is of \textit{type I} if it contains a faithful abelian projection.
A
von Neumann algebra $M$ without non zero abelian projections is
called \textit{continuous}. An arbitrary von Neumann algebra $M$ can
be decomposed in a unique way into the direct sum of von Neumann
algebras of type I$_{fin}$ (finite type I), type I$_{\infty}$
(properly infinite type I), type II$_1$ (finite continuous),
type II$_{\infty}$ (semi-finite, properly infinite, continuous)
and type III.

A linear subspace $\mathcal{D}$ in $H$ is said to be
\textit{affiliated} with $M$ (denoted as $\mathcal{D}\eta M$), if
$u(\mathcal{D})\subset \mathcal{D}$ for every unitary  $u$ in
the commutant
$$M'=\{y\in B(H):xy=yx, \,\forall x\in M\}$$ of the von Neumann algebra $M$ in $B(H)$.

A linear operator  $x$ on  $H$ with the domain  $\mathcal{D}(x)$
is said to be \textit{affiliated} with  $M$ (denoted as  $x\eta M$) if
$\mathcal{D}(x)\eta M$ and $u(x(\xi))=x(u(\xi))$
 for all  $\xi\in
\mathcal{D}(x)$ and for every unitary  $u$ in $M'$.

A linear subspace $\mathcal{D}$ in $H$ is said to be \textit{strongly
dense} in  $H$ with respect to the von Neumann algebra  $M,$ if
\begin{enumerate}
\item[1)] $\mathcal{D}\eta M;$
\item[2)] there exists a sequence of projections
$\{p_n\}_{n=1}^{\infty}$ in $P(M)$  such that
$p_n\uparrow\textbf{1},$ $p_n(H)\subset \mathcal{D}$ and
$p^{\perp}_n=\textbf{1}-p_n$ is finite in  $M$ for all
$n\in\mathbb{N}.$
\end{enumerate}

A closed linear operator  $x$ acting in the Hilbert space $H$ is
said to be \textit{measurable} with respect to the von Neumann
algebra  $M,$ if
 $x\eta M$ and $\mathcal{D}(x)$ is strongly dense in  $H.$ Denote by
 $S(M)$ the set of all measurable operators with respect to
 $M.$

A closed linear operator $x$ in  $H$  is said to
be \textit{locally measurable} with respect to the von Neumann
algebra $M,$ if $x\eta M$ and there exists a sequence
$\{z_n\}_{n=1}^{\infty}$ of central projections in $M$ such that
$z_n\uparrow\textbf{1}$ and $z_nx \in S(M)$ for all
$n\in\mathbb{N}.$

It is well-known (see e.g.~\cite{Mur}) that the set $LS(M)$ of all
locally measurable operators with respect to $M$ is a unital
*-algebra when equipped with the algebraic operations of the
strong addition and multiplication and taking the adjoint of an
operator.

 Let   $\tau$ be a faithful normal semi-finite trace on $M.$ We recall that a closed linear operator
  $x$ is said to be  $\tau$\textit{-measurable} with respect to the von Neumann algebra
   $M,$ if  $x\eta M$ and   $\mathcal{D}(x)$ is
  $\tau$-dense in  $H,$ i.e. $\mathcal{D}(x)\eta M$ and given   $\varepsilon>0$
  there exists a projection   $p\in M$ such that   $p(H)\subset\mathcal{D}(x)$
  and $\tau(p^{\perp})<\varepsilon.$
   Denote by  $S(M,\tau)$ the set of all   $\tau$-measurable operators with respect to  $M.$

The subalgebra $\mathcal{A}\subset LS(M)$ is said to be solid, if $x\in \mathcal{A},\, y\in LS(M),\,|y|\leq |x|$
implies $y\in \mathcal{A}.$

    It is well-known that $S(M)$ and $S(M, \tau)$ are solid *-subalgebras in $LS(M)$ (see~\cite{Mur}).

    Consider the topology  $t_{\tau}$ of convergence in measure or \textit{measure topology}
    on $S(M, \tau),$ which is defined by
 the following neighborhoods of zero:
$$V(\varepsilon, \delta)=\{x\in S(M, \tau): \exists e\in P(M), \tau(e^{\perp})\leq\delta, xe\in
M,  \|xe\|_{M}\leq\varepsilon\},$$  where $\varepsilon, \delta$
are positive numbers.

 It is well-known~\cite{Nel} that $S(M, \tau)$ equipped with the measure topology is a
complete metrizable topological *-algebra.

In the algebra   $S(M, \tau)$ consider the subset  $S_0(M, \tau)$
of all operators $x$ such that given any $\varepsilon>0$ there is
a projection  $p\in P(M)$ with $\tau(p^{\perp})<\infty,\,xp\in M$
and $\|xp\|_M <\varepsilon.$ Following~\cite{Str} let us call the
elements of $S_0(M, \tau)$ \textit{$\tau$-compact operators} with
respect to $M.$ It is known~\cite{Mur}, \cite{Yea} that $S_0(M,
\tau)$ is a $\ast$-subalgebra in  $S(M, \tau)$ and a bimodule over
$M,$ i.e. $ax, xa\in S_0(M, \tau)$ for all $x\in S_0(M, \tau)$ and
$a\in M.$

The following properties of the algebra $S_0(M, \tau)$  are known
(see~\cite{Bik}, \cite{Str}):

Let  $M$ be a von Neumann algebra
with a faithful normal semi-finite trace $\tau.$ Then
\begin{enumerate}
\item[1)] $S(M, \tau)=M+S_0(M, \tau);$
\item[2)] $S_0(M, \tau)$ is an ideal in $S(M, \tau).$
\end{enumerate}

Note that if the trace $\tau$ is finite then
$$S_0(M, \tau)=S(M, \tau)=S(M)=LS(M).$$

The following result describes one of the most important
properties of the algebra  $LS(M)$ (see~\cite{Mur}, \cite{Sai}).

\begin{prop}
\label{p2.1}
 Suppose that the von Neumann algebra  $M$ is the
$C^{\ast}$-product of von Neumann algebras  $M_i,$ $i\in I,$ where
$I$ is an arbitrary set of indices, i.e.
$$M=\bigoplus\limits_{i\in I}M_i=
\{\{x_i\}_{i\in I}:x_i\in M_i, i\in I, \sup\limits_{i\in
I}\|x_i\|_{M_i}<\infty\}$$ with the coordinate-wise algebraic
operations and involution and with the $C^{\ast}$-norm
$\|\{x_i\}_{i\in I}\|_{ M}=\sup\limits_{i\in I}\|x_i\|_{M_i}.$
Then the algebra  $LS(M)$ is *-isomorphic to the algebra
$\prod\limits_{i\in I}LS(M_i)$ (with the coordinate-wise
operations and involution), i.e.
$$LS(M)\cong\prod\limits_{i\in I}LS(M_i)$$
($\cong$ denotes  *-isomorphism of algebras). In particular, if $M$ is  finite, then
$$S(M)\cong\prod\limits_{i\in I}S(M_i).$$
\end{prop}

It should be noted that such isomorphisms are not valid in general for the algebras
 $S(M),$ $S(M, \tau)$   (see~\cite{Mur}).

Proposition 2.1 implies that given any family  $\{z_i\}_{i\in I}$
of mutually orthogonal central projections in $M$ with
$\bigvee\limits_{i\in I}z_i=\textbf{1}$, and a  family
of elements $\{x_i\}_{i\in I}$ in $LS(M)$, there exists a unique element $x\in LS(M)$
such that $z_i x=z_i x_i$ for all $i\in I.$ This element is
denoted by  $x=\sum\limits_{i\in I}z_i x_i.$

It is well-known (see e.g.~\cite{Seg}) that every commutative von
Neumann algebra
 $M$
is *-isomorphic to the algebra   $L^{\infty}(\Omega)=L^{\infty}(\Omega, \Sigma,
\mu)$ of all  essentially bounded measurable complex functions on a measure space
$(\Omega, \Sigma, \mu)$, and in this case  $LS(M)=S(M)\cong
L^{0}(\Omega),$   where $L^{0}(\Omega)=L^{0}(\Omega, \Sigma,
\mu)$ is the algebra of all  measurable complex functions on
$(\Omega, \Sigma, \mu).$

The following
 description of the centers of the algebras $S(M),$  $S(M, \tau)$ and  $S_0(M, \tau)$ for type I
 von Neumann algebras is very important in investigation of the structure of
 these algebras (see~\cite{Alb2}, \cite{Alb6}).

\begin{prop}\label{cccc}
Let  $M$ be a  von Neumann
algebra of type  I with  center  $Z$ and a
faithful normal semi-finite trace $\tau.$

a) If $M$ is finite, then
$Z(S(M))=S(Z)$  and  $Z(S(M,\tau))=S(Z, \tau_Z),$ where  $\tau_Z$ is the restriction of the trace
 $\tau$ on $Z;$

b) If  $M$ is of type  $I_{\infty},$ then
 the centers of the algebras $S(M)$ and  $S(M,\tau)$ coincide with $Z,$ and
 the center of the   algebra  $S_{0}(M,\tau)$ is trivial, i.e.  $Z(S_{0}(M,\tau))=\{0\}.$
\end{prop}

Let  $M$  be a von Neumann algebra of type I$_{n}$
$(n\in\mathbb{N})$ with  center  $Z.$ Then  $M$ is *-isomorphic
to the algebra  $M_n(Z)$ of $n\times n$ matrices over $Z$
(see~\cite{Sak1}, Theorem 2.3.3).

In this case the algebras $S(M, \tau)$ and $S(M)$ can be described
in the following way (see~\cite{Alb2}).

\begin{prop}
Given a  von Neumann algebra $M$
of type  $I_{n},$ $n\in\mathbb{N},$ with  a faithful normal
semi-finite trace $\tau,$  denote by $Z(S(M, \tau))$ and $Z(S(M))$  the centers
of the algebras $S(M, \tau)$ and $S(M),$
respectively.  Then $S(M, \tau)\cong M_n(Z(S(M,
\tau)))$ and  $S(M)\cong M_n(Z(S(M))).$
\end{prop}

\section{Derivations on algebras of measurable operators for type I von Neumann algebras}
\label{derI}

In this section we shall give a complete description of
derivations on the algebras $LS(M),$ $S(M),$  $S(M, \tau)$ and
$S_0(M, \tau)$ for a type I von Neumann algebra $M.$

First we shall present  results of Ber, Chilin and
Sukochev (see~\cite{Ber4, ber-mt, Ber2, Ber}) concerning the
existence of nontrivial derivations on the algebras $S(M)$ and
$S(M,\tau)$ in the case where $M$ is an abelian  von
Neumann algebra.

Let $\mathcal{A}$ be a commutative algebra with  unit $\mathbf{1}$ over
the field $\mathbb{C}$ of complex numbers. We denote by $\nabla$
the set $\{e\in \mathcal{A}: e^2=e\}$ of all idempotents in $\mathcal{A}.$ For $e,f\in
\nabla$ we set $e\leq f$ if $ef=e.$ Equipped with this partial order,
 lattice operations $e\vee f=e+f-ef, \ e\wedge f=ef$ and
the complement $e^{\bot}=\mathbf{1}-e$, the set $\nabla$ forms a
Boolean algebra. A non zero element $q$ from the Boolean algebra
$\nabla$ is called an \textit{atom} if $0\neq e\leq q, \ e\in
\nabla,$ imply that $e=q.$ If given any nonzero $e\in \nabla$
there exists an atom $q$ such that $q\leq e,$ then the Boolean
algebra $\nabla$ is said to be \textit{atomic}.

An algebra $\mathcal{A}$ is called \textit{regular} (in the sense of von Neumann)
if for any $a\in \mathcal{A}$ there exists $b \in A$ such that $a = aba.$

Along this section, we shall always assume that $\mathcal{A}$ is a
unital commutative regular algebra over $\mathbb{C},$ and that
$\nabla$ is the Boolean algebra of all its idempotents. In this
case given any element $a\in \mathcal{A}$ there exists an
idempotent $e\in \nabla$  such that $ea=a,$ and if $ga=a, g\in
\nabla,$ then $e\leq g.$ This idempotent is called the
\textit{support} of $a$ and denoted by $s(a)$ (see \cite[P.
111]{Ber}).

Suppose that $\mu$ is a strictly positive countably additive
finite measure on the Boolean algebra $\nabla$ of idempotents in
$\mathcal{A}$, and let us consider the metric
$\rho(a,b)=\mu(s(a-b)),\ a,b\in \mathcal{A},$ on the algebra  $A$.
From now on we shall assume that $(\mathcal{A}, \rho)$ is a
complete metric space (cf.~\cite{Ayu1}, \cite{Ber}).

\begin{example}
 The most important example of a complete
 commutative regular algebra $(\mathcal{A},\rho)$ is the algebra
$\mathcal{A}=L^0(\Omega)=L^0(\Omega,\Sigma,\mu)$ of all
 measurable complex functions on a measure space $(\Omega,\Sigma,\mu),$
where  $\mu$ is a finite, countably additive measure on $\Sigma,$
and $\rho(a, b)=\mu(s(a-b))=\mu(\{\omega\in \Omega:a(\omega)\neq
b(\omega)\})$ (see for details~\cite{Ayu1}, Lemma and \cite{Ber},
Example~2.2).
\end{example}

If $(\Omega,\Sigma,\mu)$ is a general
localizable measure space, i.e. the  measure $\mu$  (not finite in general)
has the finite sum property, then the algebra
$L^0(\Omega,\Sigma,\mu)$ is a unital  regular algebra, but
$\rho(a,b)=\mu(s(a-b))$
 is not a metric in general. But one can
represent $\Omega$ as a union of pair-wise disjoint measurable
sets with  finite measures and thus this algebra is a direct sum
of commutative regular complete metrizable algebras from the above example.

Following~\cite{Ber} we say that an element $a \in \mathcal{A}$ is
\textit{finitely valued} (respectively, \textit{countably valued})
if $a = \sum\limits_{k =1}^n {\alpha _k e_k}$, where $\alpha _k
\in \mathbb{C}$, $e_k \in \nabla, \ e_k e_j=0, \ k\neq j,\ k,j
=1,...,n, \ n\in \mathbb{N}$ (respectively, $a = \sum\limits_{k
=1}^\omega {\alpha _k e_k}$, where $\alpha_k \in \mathbb{C}$, $e_k
\in \nabla, \ e_k e_j=0, \ k\neq j, \ k,j =1,...,\omega,$ where
$\omega$ is a natural number or $\infty$ (in the latter case the
convergence of series is understood with respect to the metric
$\rho$)). We denote by $K(\nabla)$ (respectively, by $K_c(\nabla)$)
the set of all finitely valued (respectively, countably valued)
elements in $\mathcal{A}.$ It is known that $\nabla \subset
K(\nabla) \subset K_c(\nabla),$ and that both $K(\nabla)$ and $K_c(\nabla)$
are regular subalgebras in $\mathcal{A}.$  Moreover, the closure
of $K(\nabla)$ in $(\mathcal{A},\rho)$ coincides with
$K_c(\nabla)$ (see~\cite{Ber}, Proposition 2.8).

The following theorem provides a necessary and sufficient condition for a commutative regular algebra
to admit  nontrivial  derivations
(see~\cite{Ber2}, \cite{Ber}).

\begin{theorem}\label{bma}
$\mathcal{A}$ be a unital commutative regular algebra over $\mathbb{C}$ and let $\mu$ be a
strictly positive countably additive finite measure on the Boolean algebra $\nabla$ of all idempotents in $\mathcal{A}.$
Suppose that $\mathcal{A}$ is complete in
 the metric $\rho(a,b)=\mu(s(a-b)),\ a,b\in \mathcal{A}.$ Then the following conditions are
equivalent:
\begin{enumerate}
\item[(i)] $K_c(\nabla)\neq \mathcal{A};$
\item[(ii)] The algebra $\mathcal{A}$ admits a non-zero derivation.
\end{enumerate}
\end{theorem}

An important special case of  Theorem 3.2  is the following result
concerning the regular algebra $L^0(\Omega,\Sigma,\mu)$
(see~\cite{Ber2}, \cite{Ber}).
 \medskip

\begin{cor}\label{DIS}
Let  $(\Omega,\Sigma,\mu)$ be a finite measure space and let $L^0(\Omega)=L^0(\Omega,\Sigma,\mu)$
be the algebra of all real or complex valued measurable functions on $(\Omega,\Sigma,\mu).$
The following conditions are equivalent:
\begin{enumerate}
\item[(i)]  the Boolean algebra of all idempotents from $L^0(\Omega)$ is
not atomic;
\item[(ii)] $L^0 (\Omega)$  admits a non-zero derivation.
\end{enumerate}
\end{cor}

It is well known \cite[P. 45]{Sak1} that if $M$ is a commutative
von Neumann algebra with a faithful normal semi-finite trace
$\tau$, then $M$ is *-isomorphic to the algebra
$L^{\infty}(\Omega)=L^{\infty} (\Omega,\Sigma,\mu)$ of all
essentially bounded measurable complex valued function on an
appropriate localizable measure space $(\Omega,\Sigma,\mu)$ and
$\tau(f)=\int \limits_{\Omega}f(t)d\mu(t)$ for  $f\in L^{\infty}
(\Omega,\Sigma,\mu).$ In this case the algebra $S(M)$ of all
 measurable operators affiliated with
$M$ may be identified with the algebra $L^0(\Omega)=L^0(\Omega,\Sigma,\mu)$ of all measurable
complex valued functions on $(\Omega,\Sigma,\mu),$
while the algebra $S(M, \tau)$ of $\tau$-measurable
operators from $S(M)$ coincides with the algebra
$$\{f\in L^0 (\Omega):\exists F\in\Sigma,  \mu(\Omega \setminus F)<+\infty,  \chi_F \cdot f \in L^{\infty} (\Omega)\}$$
of all totally $\tau$-measurable functions on $\Omega,$ where $\chi_F$ is the
characteristic function of the set $F.$
If the trace $\tau$ is finite then $S(M,\tau)=S(M)\cong L^0(\Omega) $
are commutative regular algebras. But if the trace
$\tau$ is not finite, the algebra $S(M,\tau)$ is not regular. In
this case, by considering $\Omega$ as a union of pairwise disjoint
measurable sets with finite measures, we obtain that $S(M)$ is  a
direct sum of commutative regular algebras which are metrizable in the above
sense, and hence $S(M,\tau)$ is a solid subalgebra of   this direct sum.
Therefore  Corollary~\ref{DIS} implies  the following solution of
the problem concerning existence of   derivations on algebras
of measurable operator in the commutative case (see~\cite{Ber2},
\cite{Ber}).

\begin{theorem}\label{abel}
Let $M$ be a commutative von Neumann algebra with a faithful normal
semi-finite trace $\tau.$ The following conditions are equivalent:
\begin{enumerate}
\item[(i)] The lattice $P(M)$ of projections in $M$ is not atomic;
\item[(ii)]
The algebra $S(M)$ (respectively $S(M, \tau)$) admits a non-inner
derivation.
\end{enumerate}
\end{theorem}

We are now in position to give a complete description of all
derivations on the algebras $LS(M),$ $S(M),$  $S(M, \tau)$ and
$S_0(M, \tau)$ for a type I von Neumann algebra $M.$ These results
were obtained by
 Albeverio, Ayupov and Kudaybergenov (see~\cite{Alb1, Alb4, Alb2,  Alb5, Alb6},  and \cite{Ayu4}).

It is clear that if a derivation $D$ on $LS(M)$ is inner then it
is $Z$-linear, i.e. $D(f x)=f D(x)$ for all $f\in Z,$ $x\in
LS(M),$ where $Z$ is the center of the von Neumann algebra $M.$
The following main result of~\cite{Alb1} asserts that the converse
is also true.

\begin{theorem}\label{tttt}
Let  $M$ be a type I von Neumann algebra with center $Z.$ A derivation
 $D$ on the algebra $LS(M)$ is inner if and only if it is
 $Z$-linear, or equivalently it is identically zero on $Z.$
 \end{theorem}

Let $\mathcal{A}$ be a commutative algebra and let $M_n(\mathcal{A})$
 be the algebra of $n\times n$ matrices over $\mathcal{A}.$
If  $e_{i,j},\,i,j={1,..., n},$
are the matrix units in  $M_n(\mathcal{A}),$ then each element
$x\in M_n(\mathcal{A})$ has the form
$$
x=\sum\limits_{i,j=1}^{n}f_{i j}e_{i j},\,f_{i,j}\in \mathcal{A},\,i,j=1,2,..., n.
$$
Let  $\delta:\mathcal{A}\rightarrow \mathcal{A}$ be a
derivation. Setting
\begin{equation}
\label{1} D_{\delta}\left(\sum\limits_{i,j=1}^{n}f_{i j}e_{i
j}\right)= \sum\limits_{i,j=1}^{n}\delta(f_{i j})e_{i j}
\end{equation}
we obtain a well-defined linear operator $D_\delta$ on the algebra
$M_n(\mathcal{A}).$ Moreover $D_\delta$ is a derivation on the
algebra  $M_n(\mathcal{A})$, and its restriction onto the center of
the algebra  $M_n(\mathcal{A})$ coincides with the given $\delta.$

 Now let us consider arbitrary (non $Z$-linear, in general)
derivations on $LS(M)$. The following simple but important remark is crucial in our further considerations.

Let  $\mathcal{A}$ be an algebra with  center  $Z$ and let
$D:\mathcal{A}\rightarrow \mathcal{A}$ be a derivation. Given any
$x\in \mathcal{A}$ and a central element $f \in Z$  we have
$$
D(f x)=D(f)x+f D(x)
$$
and
$$
D(xf)=D(x)f +xD(f).
$$
Since  $f x=xf$ and $f D(x)=D(x)f,$ it follows that $D(f)x=xD(f)$
for any $x\in \mathcal{A}.$ This means that $D(f)\in Z,$ i.e.
$D(Z)\subseteq Z.$ Therefore, given any derivation $D$ on the
algebra $A$ we can consider its restriction $\delta:Z\rightarrow
Z.$

Now let  $M$ be a homogeneous von Neumann algebra of type I$_{n},
n \in \mathbb{N}$,  with center $Z.$ Then  the algebra $M$ is *-isomorphic
to the algebra $M_n(Z)$ of all  $n\times n$- matrices over $Z,$
and the algebra  $LS(M)=S(M)$ is *-isomorphic to the algebra
 $M_n(S(Z))$ of all  $n\times n$ matrices over
$S(Z),$ where $S(Z)$ is the algebra of measurable operators with respect to the commutative von Neumann algebra $Z$.

The algebra $LS(Z)=S(Z)$ is isomorphic to the algebra $L^{0}(\Omega)=L^{}(\Omega, \Sigma, \mu)$
 of all measurable complex functions on a measure space,
 and therefore it admits (in non atomic cases) non zero derivations
  (see Theorem~\ref{abel}).

The following consideration is the main step in constructing
 the ''exotic'' derivation $D_\delta$
on the algebra $S(M)$ of measurable operators
affiliated with a finite type I von Neumann algebra $M$, which
admits a non trivial derivation $\delta$ on its center $S(Z)$.

Let  $\delta:S(Z)\rightarrow S(Z)$ be a
derivation and  let $D_\delta$ be the derivation on the algebra $M_n(S(Z))$ defined by (\ref{1}).

The following lemma describes the structure of an arbitrary
derivation on the algebra of locally measurable operators for
homogeneous type  I$_n,$ $n\in\mathbb{N},$ von Neumann algebras
(see~\cite{Alb2}).

\begin{lemma}\label{ddddd}
Let  $M$ be a homogenous von Neumann algebra of type
  $I_{n}, n \in \mathbb{N}.$
Every derivation  $D$ on the algebra $LS(M)$ can be uniquely represented as a sum
  $$D=D_{a}+D_{\delta ,}$$ where  $D_{a}$ is an inner derivation implemented by an element  $a\in LS(M),$
while $D_{\delta} $ is the derivation of the form (\ref{1}), generated by
a derivation $\delta$ on the center of $LS(M)$ identified
with $S(Z)$.
\end{lemma}

Now let  $M$ be an arbitrary finite von Neumann algebra of type I
with center $Z.$ There exists a family  $\{z_n\}_{n\in F},$
$F\subseteq\mathbb{N},$ of central projections from $M$ with
$\sup\limits_{n\in F}z_n=\textbf{1}$, such that the algebra  $M$ is
*-isomorphic to the  $C^{*}$-product of von Neumann algebras
$z_n M$, where each $z_n M$ is of type  I$_{n}$ respectively, $n\in F,$ i.e.
$$M\cong\bigoplus\limits_{n\in F}z_n M.$$
By Proposition 2.1 we have that
$$LS(M)\cong\prod\limits_{n\in F}LS(z_n M).$$

Suppose that    $D$ is a derivation on  $LS(M),$ and
$\delta$ is its restriction onto its center  $S(Z).$
Since  $\delta$ maps each  $z_nS(Z)\cong Z(LS(z_n M))$ into itself, for each $n$,  it generates a derivation $\delta_n$
on $z_nS(Z)$ for each $n\in F.$

Let     $D_{\delta_n}$ be the derivation on the matrix algebra
$M_n(z_nZ(LS(M)))\cong LS(z_nM)$ defined as in
(\ref{1}). Put
\begin{equation}
\label{2}
D_\delta(\{x_n\}_{n\in F})=\{D_{\delta_n}(x_n)\},\,\{x_n\}_{n\in F}\in LS(M).
\end{equation}
Then the map  $D$  is a derivation on $LS(M).$

Now Lemma~\ref{ddddd} implies the following result, which shows,
in particular, that  $D_\delta$ is the most general
form of non-inner derivations on $LS(M).$

\begin{lemma}
Let  $M$ be a finite von Neumann algebra
of type I. Each
derivation  $D$ on the algebra  $LS(M)$ can be uniquely
represented in the form
$$D=D_{a}+D_{\delta ,}$$
where  $D_{a}$ is an inner derivation implemented by an element
$a\in LS(M),$ and $D_{\delta} $ is a derivation given as in
(\ref{2}).
\end{lemma}

Now we shall  consider derivations on  algebras of locally measurable
operators affiliated with type I$_{\infty}$ von Neumann algebras.

\begin{theorem}\label{typei}
If $M$ is a type $I_{\infty}$ von Neumann algebra, then any
derivation on the algebras $LS(M),$ $S(M)$ and $S(M, \tau)$  is inner.
\end{theorem}

Finally, let us consider derivations on the
algebra $LS(M)$ of locally measurable operators with respect to
an arbitrary  type I von Neumann algebra $M.$

Let $M$ be a type  I von Neumann algebra. There exists a central
projection  $z_0\in M$ such that
\begin{enumerate}
\item[a)] $z_0M$ is a finite von Neumann algebra;
\item[b)] $z_0^{\bot}M$ is a von Neumann algebra of type  I$_{\infty}.$
\end{enumerate}

Consider a derivation  $D$ on  $LS(M)$ and let $\delta$ be its
restriction onto its center  $Z(S).$ By Theorem~\ref{typei} the restriction
 $z_0^{\bot}D$ of the derivation $D$ onto $z_0^{\bot}LS(M)$ is inner,
  and thus  we have $z_0^{\bot}\delta\equiv 0,$ i.e. $\delta=z_0\delta.$

Let   $D_\delta$ be the derivation on  $z_0LS(M)$ defined as in
\eqref{2} and consider its extension $D_\delta$ on
$LS(M)=z_0LS(M)\oplus z_0^{\bot}LS(M)$, which is defined as
\begin{equation}
\label{3}
D_\delta(x_1+x_2):=D_\delta(x_1),\,x_1\in z_0LS(M),x_2\in z_0^{\bot}LS(M).
\end{equation}

The following theorem is the main result of this section, and
gives the general form of derivations on the algebra $LS(M)$
(see~\cite{Alb2}).

\begin{theorem}\label{typee}
Let  $M$ be a type  $I$ von Neumann
algebra and let $A$ be one of the algebras $LS(M),$  $S(M)$ or $S(M, \tau).$ Each
derivation  $D$ on  $A$ can be uniquely represented in
the form
\begin{equation}\label{D}
   D=D_{a}+D_{\delta}
\end{equation}
where  $D_{a}$ is an inner derivation implemented by an element
$a\in A,$ and $D_{\delta} $ is a derivation of the form
(\ref{3}), generated by a derivation $\delta$ on the center of $A.$
\end{theorem}

 If we consider the measure topology $t_{\tau}$ on the algebra $S(M, \tau)$
 then it is clear that every non-zero derivation of the form
  $D_\delta$ is discontinuous in $t_\tau.$ Therefore the above Theorem~\ref{typee} implies:

\begin{cor}\label{HOM}
Let $M$ be a type  I von Neumann algebra with a faithful normal semi-finite trace $\tau.$ A derivation $D$
on the algebra    $S(M, \tau)$ is inner if and only if it is continuous in the measure topology.
\end{cor}

Now, let $M$ be a type I von Neumann algebra with atomic
center $Z$ and let $\{q_i\}_{i\in I}$ be the set of all atoms of
$Z.$ Consider a derivation $D$ on $LS(M).$ Since $q_iZ\cong
q_i\mathbb{C}$ for all $i\in  I,$ we have $q_iD(f x) = D(q_i f x)
= q_if D(x)$ for all $i\in  I, f\in Z, x\in LS(M).$ Thus $D(f x) =
f D(x)$ for all $f\in Z.$ This means that in the case of $Z$  being atomic,
 every derivation on $LS(M)$
 is automatically $Z$-linear. Combining this fact with
Theorem~\ref{tttt}, we have the following result which is a
strengthening of  result obtained by Weigt in~\cite{Wei}.

\begin{cor}\label{atle}
If $M$ is a von Neumann algebra with  atomic lattice of projections, then
every derivation on the algebras  $LS(M), S(M)$ and $S(M, \tau)$ is inner.
\end{cor}

Now let us consider derivations on the algebra $S_0(M, \tau)$ of
$\tau$-compact operators affiliated with a semi-finite von
Neumann algebra $M$ and a faithful normal semi-finite trace $\tau$
(see~\cite{Alb4}, \cite{Alb6}).

It should be noted that for an arbitrary von Neumann algebra  $M,$
the center of the algebra $LS(M)$  coincides with $LS(Z)$, and thus
contains $Z$ (see Proposition~\ref{cccc}). This was an essential
point in the proof of theorems  describing
derivations on the algebra $LS(M)$ of locally measurable operators
with respect to a type I von Neumann algebra $M.$ Proposition 2.2
shows that this is not the case for the algebra  $S_{0}(M,\tau)$,
because the center of this algebra may be trivial. Thus, the
methods of the proof of Theorem 4.1 from~\cite{Alb2} can not be
directly applied for  description of derivations on algebras
of $\tau$-compact operators with respect to type I von Neumann
algebras. Nevertheless, the following result for the algebra
$S_0(M, \tau)$ is obtained in~\cite{Alb6}.

\begin{theorem}
Let  $M$ be a type  I von Neumann
algebra with a faithful normal semi-finite trace $\tau.$ Each
derivation  $D$ on  $S_0(M, \tau)$ can be uniquely represented in
the form
$$D=D_{a}+D_{\delta},$$
where  $D_{a}$ is a spatial derivation implemented by an element
$a\in S(M, \tau),$ and $D_{\delta} $ is a derivation of the form
(\ref{3}), generated by a derivation $\delta$ on the center of
$S_0(M, \tau)$.
\end{theorem}

Recently, in \cite{AKIE13} we have investigated derivations on
algebras of $\tau$-compact operators affiliated with an arbitrary
 semi-finite (i.e. type II algebras are also included) von
Neumann algebra $M$ and a faithful normal semi-finite trace $\tau$.
Namely, we proved that every $t_\tau$-continuous
derivation on the algebra $S_0(M, \tau)$ is spatial and
implemented by a $\tau$-measurable operator affiliated with $M,$
where $t_\tau$ denotes the measure topology on  $S_0(M, \tau)$. We
have also shown  automatic $t_\tau$-continuity of all
derivations on $S_0(M, \tau)$ for properly infinite von Neumann
algebras $M$. Thus, in the properly infinite case the condition of
$t_\tau$-continuity of the derivation  is redundant for its
spatiality.

\begin{theorem}\label{ieot}
Let $M$ be a    von Neumann algebra with a faithful normal
semi-finite trace $\tau.$ Then every $t_\tau$-continuous
derivation $D:S_0(M, \tau)\rightarrow S_0(M, \tau)$ is spatial and
implemented by an element $a\in S(M, \tau).$
\end{theorem}

\begin{theorem}\label{contin}
 Let  $M$ be a properly infinite von Neumann algebra with a faithful normal
semi-finite trace $\tau.$ Then any derivation $D:S_0(M,
\tau)\rightarrow S_0(M, \tau)$  is $t_\tau$-continuous.
\end{theorem}

From Theorems~\ref{ieot} and  \ref{contin}   we obtain the
following result.

\begin{theorem}\label{proper}
 If $M$ is a properly infinite von Neumann algebra with a faithful normal
semi-finite trace $\tau,$ then any derivation $D:S_0(M,
\tau)\rightarrow S_0(M, \tau)$  is  spatial and implemented by an
element $a\in S(M, \tau).$
\end{theorem}

\section{Derivations on algebras of measurable operators for arbitrary von Neumann algebras}
\label{derG}

In the present section  we shall consider derivations on the
algebras $LS(M)$ and $S(M)$ for an arbitrary von Neumann algebra
$M.$

First we consider additive derivations on the algebra $LS(M),$
where $M$ is a properly infinite von Neumann algebra. These
results are obtained in the paper of Ayupov and
Kudaybergenov (see~\cite{Ayu2, AK2013}).

We  shall consider  the so called  central extension $E(M)$ of a
von Neumann algebra $M$ and   show that
 $E(M)$ is a *-subalgebra in the algebra $LS(M)$ and
 this subalgebra coincides with whole $LS(M)$ if and only if
 $M$ does not contain a direct summand of type II.

 As the main result of this section we obtain
  that if  $M$ is a properly infinite  von
Neumann algebra, then every additive derivation on the algebra
$E(M)$ is inner. In particular, every additive derivation on the
algebra $LS(M),$ where $M$ is of type I$_\infty$ or III, is inner.

Let  $E(M)$   denote the set of all elements  $x$ from  $LS(M)$
for which there exists a sequence of mutually orthogonal central
projections  $\{z_i\}_{i\in I}$ in  $M$ with $\bigvee\limits_{i\in
I}z_i=\textbf{1},$ such that $z_i x\in M$ for all $i\in I,$ i.e.
$$
E(M)=\{x\in LS(M): \exists\, z_i\in P(Z(M)), z_iz_j=0, i\neq j,
\bigvee\limits_{i\in I}z_i=\textbf{1}, z_i x\in M, i\in I\},
$$
where $Z(M)$ is the center of $M.$

\begin{prop}\label{cee}
Let $M$ be a von Neumann algebra with the center $Z(M).$
 Then
\begin{enumerate}
\item[i)]  $E(M)$ is a *-subalgebra in $LS(M)$ with  center
$S(Z(M)),$ where $S(Z(M))$ is the algebra of measurable operators
with respect to  $Z(M);$
\item[ii)] $LS(M)=E(M)$ if and only if $M$ does not have
direct summands of type II.
\end{enumerate}
\end{prop}

A similar notion (i.e. the algebra $E(\mathcal{A})$) for arbitrary
*-subalgebras $\mathcal{A}\subset LS(M)$  was independently
introduced recently by Muratov and Chilin
in~\cite{Mur1}. They called it the central extension of
$\mathcal{A}.$ Therefore following~\cite{Mur1} we shall say that
$E(M)$ is \textit{the central extension of} $M$.

One has the following description of $E(M)$ (see~\cite{Ayu2},
\cite{Mur1}).

\begin{prop}\label{nnn}
Let $M$ be a von Neumann algebra. Then $x\in E(M)$ if only if
there exists  $f\in S(Z(M))$ such that  $|x|\leq f.$
\end{prop}

The following theorem is obtained in \cite{Ayu2}).

\begin{theorem}\label{mmm}
Let  $M$ be a properly infinite  von
Neumann algebra. Then
every additive derivation on the algebra $E(M)$ is inner.
\end{theorem}

The proof of Theorem~\ref{mmm} is based on the following lemma which  has
some interest in its own right
(see~\cite{Ayu2}, \cite{Ber3}).

\begin{lemma}
Let   $M$ be a properly infinite von
Neumann algebra, and let $\mathcal{A}\subseteq LS(M)$ be a *-subalgebra such
that  $M\subseteq \mathcal{A}$, and suppose that $D:\mathcal{A}\rightarrow \mathcal{A}$ is an
additive derivation. Then  $D|_{Z(\mathcal{A})}\equiv 0,$ in particular,
$D$ is $Z(\mathcal{A})$-linear.
\end{lemma}

From Theorem~\ref{mmm}  and Proposition \ref{cee} we obtain the
following extension of Theorem~\ref{typei}.

\begin{theorem}\label{mixx}
Let  $M$ be a  direct sum of von Neumann algebras
of type I$_\infty$ and III. Then every additive derivation on
the algebra $LS(M)$ is inner.
\end{theorem}

Since $LS(M)$ contains $S(M)$ as a solid *-subalgebra, and
$S(M)$ contains $S(M, \tau)$ as a solid *-subalgebra,  Theorem~\ref{mixx} implies
 similar results for derivations on the algebras $S(M)$ and $S(M,\tau)$ for
type I and type III von Neumann algebras  $M.$

Thus,  the problem of describing the
derivations on the above algebras is reduced to the case, where the underlying  von
Neumann algebra is of type II.

Recently, Ber, Chilin and Sukochev in
\cite{Ber2013} have proved that any derivation on the algebra
$LS(M)$ of all locally measurable operators affiliated with a
properly infinite von Neumann algebra  $M$     is continuous with
respect to  so-called local measure topology. For type I and
type III cases this follows from our Theorem~\ref{mixx}.   But this
result is new for the type II$_\infty$ case.
Later in \cite{Ber2014} they proved the following extension of our
Theorem~\ref{mixx} for the type     II$_\infty$ case.

\begin{theorem}\label{bcs}  Every derivation on the algebra $LS(M)$ is
inner, provided that $M$ is a properly infinite von Neumann
algebra.
\end{theorem}

Therefore, the problem remains unsolved only in the case when $M$ is  a type II$_1$   von Neumann
algebra. A partial answer for this case is given by the following theorem from \cite{AK2013}.

\begin{theorem}\label{ak13}
 Let $M$ be a finite von Neumann algebra with a faithful normal
semi-finite trace  $\tau,$ equipped with the local measure
topology $t.$ Then every $t$-continuous derivation
$D:S(M)\rightarrow S(M)$ is inner.
\end{theorem}

The above theorem follows also from the above mentioned paper of
Ber, Chilin and Sukochev  in \cite{Ber2014}.

Thus the  problem of  innerness of
derivations on algebras of measurable operators is open only for the case of
type II$_1$   von Neumann algebras. For finite von Neumann algebras, the above algebras
 $S(M,\tau),$ $S(M),$ $LS(M)$ coincide with the algebra of all
closed operators affiliated with  $M$ (this is so called
Murray--von Neumann algebra) (see also \cite{KL}).

\begin{prob}
Let $M$ be a type II$_1$ von Neumann algebra (in particular -- a
II$_1$-factor). Prove that  every  derivation on the
algebra $LS(M)=S(M)$ is inner, or give an example of a $t$-discontinuous
derivation on $S(M)$ .
\end{prob}

\section{Local derivations on algebras of measurable operators}
\label{derloc}

In this section  we study  local derivations on the algebra
$S(M,\tau)$ of $\tau$-measurable operators affiliated with a von
Neumann algebra $M$ and a faithful normal semi-finite trace
$\tau.$ The results presented here are due to Albeverio, Ayupov,
Kudaybergenov and Nurjanov (see~\cite{Nur},
\cite{Nur1}).

There exist various types of linear operators which are close to
derivations (see e.g.~\cite{Bre3, Bre1, Kad, Lar}). In particular
Kadison  introduced and investigated in~\cite{Kad}   so-called
local derivations on Banach algebras.

A linear operator $\Delta$ on an algebra $A$ is called
\textit{local derivation} if given any $x\in A$ there exists a
derivation $D$ (depending on $x$) such that $\Delta(x)=D(x).$  The
main problem concerning these operators is to find conditions under
which local derivations become derivations and to present examples
of algebras which admit local derivations that are not derivations (see
e.g.~\cite{Kad}, \cite{Lar}). In particular Kadison in~\cite{Kad}
proves that every continuous local derivation from a von
Neumann algebra $M$ into a dual $M$-bimodule is a derivation.
Later this result has been extended in~\cite{Bre3} to a larger
class of linear operators $\Delta$ from $M$ into a normed
$M$-bimodule $E$ satisfying the identity
\begin{equation}
\label{5}
\Delta(p)=\Delta(p)p + p\Delta(p)
\end{equation}
for every idempotent $p\in M.$

It is clear that each local derivation satisfies (\ref{5}) since given any idempotent $p\in M$, we have
$\Delta(p)=D(p)=D(p^2)=D(p)p + pD(p)= \Delta(p)p + p\Delta(p).$

In~\cite{Bre1}  {Bre\v{s}ar and \v{S}emrl} proved that every linear operator $\Delta$
on the algebra $M_{n}(R)$ satisfying (\ref{5}) is automatically a
derivation, where $M_{n}(R)$ is the algebra of $n\times n$
matrices over a unital commutative ring $R$ containing $1/2.$

In~\cite{Joh} Johnson extends Kadison's result and proves every
local derivation from a $C^{\ast}$-algebra $\mathcal{A}$ into any
Banach $\mathcal{A}$-bimodule is a derivation. He also shows that
every local derivation from
 a $C^{\ast}$-algebra $\mathcal{A}$ into any Banach $\mathcal{A}$-bimodule is bounded
 (see~\cite[Theorem 5.3]{Joh}).

In~\cite{Tim} it was proved that every local derivation on the
maximal $O^{\ast}$-algebra $\mathcal{L}^{+}(\mathcal{D})$ is an
inner derivation.

In the present section  we study local derivations
on the algebra $S(M,\tau)$ of all $\tau$-measurable operators
affiliated with a von Neumann algebra $M$ and a faithful normal
semi-finite trace $\tau.$ One of our main results
(Theorem~\ref{seven}) presents an unbounded version of Kadison's
Theorem A from~\cite{Kad}, and it asserts that every local
derivation on $S(M,\tau)$ which is continuous in the measure
topology automatically becomes a derivation. In particular,
for  type I von Neumann algebras $M$ all such local
derivations on $S(M,\tau)$ are inner derivations. We also show
that for type I finite von Neumann algebras without abelian direct
summands, as well as for von Neumann algebras with the atomic
lattice of projections, the continuity condition on local
derivations in the above theorem  is redundant.

Finally,   we consider the problem of existence of  local
derivations which are not derivations on algebras of measurable
operators. The consideration of such examples on various  finite-
and infinite dimensional algebras was initiated  by Kadison,
Kaplansky and Jensen (see~\cite{Kad}). We consider this problem on
a class of commutative regular algebras, which includes the
algebras of measurable functions on a finite measure space, and
obtain necessary and sufficient conditions  for the algebras of
measurable and $\tau$-measurable operators affiliated with a
commutative von Neumann algebra to admit local derivations which
are not derivations.

Recall  that $S(M,\tau)$ is a complete metrizable topological $*$-algebra with
respect to the measure topology $t_{\tau}.$ Moreover, the
 algebra $S(M,\tau)$ is \textit{semi-prime}, i.e.
$aS(M,\tau)a=\{0\}$ for $a\in S(M,\tau),$ implies $a=0.$

One of  the main results of this section is the following
(see~\cite{Nur}).

\begin{theorem}\label{seven}
Let $M$ be a von Neumann algebra with a
faithful normal semi-finite trace $\tau.$ Then every $t_{\tau}$-continuous linear operator
$\Delta$ on the algebra $S(M,\tau)$ satisfying the identity (\ref{5})
is a derivation on $S(M,\tau).$ In particular any $t_{\tau}$-continuous local
derivation on the algebra $S(M,\tau)$ is a derivation.
\end{theorem}

It should be noted that the proof of the latter theorem
essentially relies on a result of  {Bre\v{s}ar }\cite[Theorem 1]{Bre2} which
asserts that every Jordan derivation on a semi-prime algebra is a
(associative) derivation.

For type I von Neumann algebras the above result can be strengthened as follows.

\begin{cor}
Let $M$ be a type I von Neumann algebra with
a faithful normal semi-finite trace $\tau.$ Then every $t_{\tau}$-continuous linear operator $\Delta$
satisfying (\ref{5}) (in particular every $t_{\tau}$-continuous  local derivation)
on $S(M,\tau)$ is an inner derivation.
\end{cor}

Further we have the following technical result, which has some intrinsic interest.
 \medskip

\begin{lemma}\label{zzz}
Every local derivation $\Delta$ on the algebra $S(M,\tau)$
is necessarily $P(Z)$-homogeneous, i.e. $$\Delta(zx)=z\Delta(x)$$
for any central projections $z\in P(Z)=P(M)\cap Z,$ and for all $x\in S(M,\tau).$
\end{lemma}

For finite von Neumann algebras the condition of $t_{\tau}$-continuity of
local derivations can be omitted. Namely, one has the following theorem.

\begin{theorem}\label{ldf}
Let $M$ be a finite von
Neumann algebra of type I without abelian direct summands, and let
$\tau$ be a faithful normal semi-finite trace on $M.$
Then every local derivation $\Delta$ on the algebra
$S(M,\tau)$ is a derivation, and hence can be represented
as the sum  (\ref{D}) of an inner derivation and a discontinious derivation.
\end{theorem}

Recently similar problems in a more general setting were also
considered  by Hadwin and coauthors in \cite{Had}. In particular, Theorem 1 from
\cite{Had} implies the following extension of the above theorems for general von Neumann algebras.

\begin{theorem}\label{haqwin}
Let $M$ be a  von Neumann algebra without abelian direct summands,
and let $\mathcal{A}$ be a subalgebra in $LS(M)$ such that $M\subseteq
\mathcal{A}.$ Then every local derivation $\Delta$ on
$\mathcal{A}$ is a derivation.
\end{theorem}

  In the latter theorems the condition on $M$ to have no abelian direct
summand is crucial, because in the case of abelian von Neumann algebras the picture is completely different.
Therefore, below  we shall consider local derivations on the algebras of measurable and $\tau$-measurable
operators affiliated with abelian von Neumann algebras.

Now let $D$ be a derivation on a regular commutative algebra
$\mathcal{A}.$ Since any derivation on $\mathcal{A}$ does  not
enlarge the supports of elements (see~\cite[Theorem]{Ayu1} and
~\cite[Proposition 2.3]{Ber})  we have that $s(D(a))\leq s(a)$ for
any $a\in \mathcal{A},$ and also $D|_{\nabla}=0.$ Therefore by the
definition, each local derivation $\Delta$ on $\mathcal{A}$
satisfies the following two conditions:
\begin{equation}
\label{A}
s(\Delta(a))\leq s(a), \  \forall\, a\in \mathcal{A},
\end{equation}
\begin{equation}
\label{B}
\Delta|_{\nabla}\equiv0.
\end{equation}
This means that~\eqref{A}  and~\eqref{B}  are necessary conditions for a linear operator
 $\Delta$ to be a local derivation on the algebra $\mathcal{A}.$

The following lemma which assert
 that  these two condition are also sufficient, is the crucial step for the proofs of the further  results in this section.

\begin{lemma}
Each linear operator on the algebra $\mathcal{A}$ satisfying the conditions~\eqref{A}  and~\eqref{B}
is a local derivation on $\mathcal{A}.$
\end{lemma}

The following theorem presents conditions for  existence of
local derivations that are not derivations on commutative regular
algebras (cf. Theorem~\ref{bma}).
 \medskip

\begin{theorem}
Let $\mathcal{A}$ be a unital commutative regular algebra over $\mathbb{C}$, and let $\mu$ be a
strictly positive countably additive finite measure on the Boolean algebra $\nabla$ of all idempotents in $\mathcal{A}.$
Suppose that $\mathcal{A}$ is complete with respect to the metric $\rho(a,b)=\mu(s(a-b)),\ a,b\in \mathcal{A}.$
 Then the following conditions are
equivalent:
\begin{enumerate}
\item[i)] $K_c(\nabla)\neq \mathcal{A};$
\item[(ii)] The algebra $\mathcal{A}$ admits a non-zero derivation;
\item[(iii)] The algebra $\mathcal{A}$ admits a non-zero local derivation;
\item[(iv)] The algebra  $\mathcal{A}$ admits a local derivation which is not a derivation.
\end{enumerate}
\end{theorem}

 The proof of the above theorem is based on the following
 technical  result, which is the main tool for construction of local derivations which are not derivations.

\begin{lemma}
If $D$ is a derivation on a commutative regular algebra $\mathcal{A},$ then $D^2$ is a
derivation if and only if
$D\equiv 0.$
\end{lemma}

An important special case of the latter theorem is the following
result concerning the regular algebra $L^0(\Omega,\Sigma,\mu).$

\begin{cor}
Let  $(\Omega,\Sigma,\mu)$ be a finite measure space and let $L^0(\Omega)=L^0(\Omega,\Sigma,\mu)$
be the algebra of all real or complex valued measurable functions on $(\Omega,\Sigma,\mu).$
The following conditions are equivalent:
\begin{enumerate}
\item[(i)]  The Boolean algebra of all idempotents from $L^0(\Omega)$ is
not atomic;
\item[(ii)] $L^0 (\Omega)$  admits a non-zero derivation;
\item[(iii)] $L^0 (\Omega)$  admits a non-zero local derivation;
\item[(iv)] $L^0 (\Omega)$  admits a local derivation which is not a derivation.
\end{enumerate}
\end{cor}

For general commutative von Neumann algebras one has the following
result (cf. Theorem~\ref{abel}).

\begin{theorem}
Let $M$ be a commutative von Neumann algebra with a faithful normal
semi-finite trace $\tau.$ The following conditions are equivalent:
\begin{enumerate}
\item[(i)] The lattice $P(M)$ of projections in $M$ is not atomic;
\item[(ii)]
The algebra $S(M)$ (respectively $S(M, \tau)$) admits a non-inner derivation;
\item[(iii)]  The algebra $S(M)$ (respectively
$S(M, \tau)$) admits a non-zero local derivation;
\item[(iv)] The algebra $S(M)$ (respectively $S(M, \tau)$)
admits a local derivation which is not a derivation.
\end{enumerate}
\end{theorem}

\section{$2$-Local derivations on algebras of measurable operators}
\label{der2loc}

This section  is devoted to  $2$-local derivations on the algebra
$S(M)$ of measurable operators affiliated with a von Neumann
algebra $M$ of type I.  The  results presented here are due to Ayupov,
Kudaybergenov and Alauadinov (see~\cite{AKAA,AKA,AKNA}).

In 1997, Semrl \cite{Sem1}  introduced the concepts of
$2$-local derivations and $2$-local automorphisms. A  map
$\Delta:\mathcal{A}\rightarrow\mathcal{A}$  (not linear in
general) is called a
 $2$-\emph{local derivation} if  for every $x, y\in \mathcal{A},$  there exists
 a derivation $D_{x, y}:\mathcal{A}\rightarrow\mathcal{A}$
such that $\Delta(x)=D_{x, y}(x)$  and $\Delta(y)=D_{x, y}(y).$
In this paper he described $2$-local derivations and
automorphisms of the algebra $B(H)$ of
 all bounded linear operators on the infinite-dimensional
separable Hilbert space $H.$  A similar description for the
finite-dimensional case appeared later in \cite{Kim}. In our paper
\cite{JMAA} we have considered $2$-local derivations on the
algebra $B(H)$ of all linear bounded operators on an arbitrary (no
separability is assumed) Hilbert space $H$ and proved that every
$2$-local derivation on $B(H)$ is a derivation. Recently, we have extended this result
for arbitrary von Neumann algebras \cite{POST}.
 Zhang and Li \cite{Zhang} described $2$-local  derivations on
symmetric digraph algebras  and constructed  a $2$-local derivation
 which is not a derivation on  the algebra
of all upper triangular complex $2\times 2$-matrices.

Throughout this section   $\mathcal{A}$ is a unital commutative regular
algebra over $\mathbf{C},$  $\nabla$ is the Boolean algebra of all
its idempotents and  $\mu$ is a strictly positive countably
additive finite measure on  $\nabla$.  Consider the metric
$\rho(a,b)=\mu(s(a-b)),\ a,b\in \mathcal{A},$ on the algebra
$\mathcal{A}$,  and assume that
$(\mathcal{A},\rho)$ is a complete metric space (cf. \cite{Ber}).

The following  Theorem (see \cite[Theorem 3.5]{AKA}) gives a
solution of the problem concerning  existence of $2$-local
derivations which are not derivations on algebras of measurable
operator in the abelian case.

\begin{theorem}\label{vonab}
Let $M$ be an abelian von Neumann algebra. The following
conditions are equivalent:
\begin{enumerate}
\item[(i)] The lattice $P(M)$ of projections in $M$ is not atomic;
\item[(ii)] The algebra $S(M)$ admits a $2$-local derivation which is not
a derivation.
\end{enumerate}
\end{theorem}

Further in this section we shall investigate $2$-local derivations
on matrix algebras over commutative regular algebras.

Let  $M_n(\mathcal{A})$
 be the algebra of $n\times n$ matrices over a commutative regular algebra~$\mathcal{A}.$
 The following result from \cite{AKA} shows that for $n\geq2$ this algebra has a completely different property
 compared with the corresponding property of the algebra $\mathcal{A}$ in the previous Theorem.

\begin{theorem}\label{Main}
Every $2$-local derivation $\Delta: M_n(\mathcal{A})\rightarrow
M_n(\mathcal{A}),$ $n\geq2$,   is a derivation.
\end{theorem}

The proof of Theorem \ref{Main} consists of several Lemmata.

For $x\in M_n(\mathcal{A})$ by $x_{ij}$ we denote the $(i,
j)$-entry of $x,$ i.e. $e_{ii}xe_{jj}=x_{ij}e_{ij},$ where $1\leq
i,j\leq n.$

\begin{lemma}\label{AA}   For every
  $2$-local derivation $\Delta$ on  $M_n(\mathcal{A}),$ $n\geq2,$
  there exists a derivation  $D$
  such that $\Delta(e_{i j})=D(e_{i j})$ for all $i,  j\in {1,2,..., n}.$
\end{lemma}

\begin{lemma}\label{DER}
If $\Delta(e_{ij})=0$ for all $i,  j\in {1,2,..., n},$  then the
restriction $\Delta|_{\mathcal{A}}$  is a derivation.
\end{lemma}

\begin{lemma}\label{ma}
If $\Delta|_{\mathcal{A}}\equiv 0$ and $\Delta(e_{ij})=0$ for all
 $i,  j\in {1,2,..., n},$  then
  $\Delta\equiv 0.$
\end{lemma}

Now we outline the  sketch of the proof for this
Theorem~\ref{Main}.

First,  according to Lemma~\ref{AA}, one can find
a derivation $D$ on $M_n(\mathcal{A})$ such that
$(\Delta-D)(e_{ij})=0$ for all $i,  j\in {1,2,..., n}.$
Further, by Lemma \ref{DER} $\delta=(\Delta-D)|_{\mathcal{A}}$ is
a derivation.  Finally, passing to  the $2$-local derivation
$\Delta_0=\Delta-D-D_\delta$ and taking into account that
$\Delta_0(e_{ij})=0$ for all  $i,  j\in {1,2,..., n},$ and that
$\Delta_0|_{\mathcal{A}}=0,$ by Lemma~\ref{ma} we obtain that
$\Delta_0=0,$ i.e. $\Delta=D+D_\delta$ is a derivation.

Let $M$ be a von Neumann algebra and denote by $S(M)$ the algebra
of all measurable operators and by $LS(M)$ -- the algebra of all
locally measurable operators affiliated with $M.$
Theorem~\ref{Main} implies the following result.

\begin{theorem}
Let  $M$ be a finite  von Neumann algebra of type I without
abelian direct summands. Then every  $2$-local derivation on the
algebra $LS(M)=S(M)$ is a derivation.
\end{theorem}

\begin{theorem}\label{Final}
Let $M$ be an arbitrary von Neumann algebra of type I$_\infty$ and
let $\mathcal{B}$ be a *-subalgebra of $LS(M)$ such that
$M\subseteq \mathcal{B}.$
 Then every $2$-local
derivation $\Delta: \mathcal{B}\rightarrow \mathcal{B}$
 is a
derivation.
\end{theorem}

The proof of Theorem~\ref{Final}  (see ~\cite{AKAA}) is
essentially different compared with the proof in the case of
finite type I von Neumann algebras. In this case we use the
extended center valued trace  $\Phi$ on the set $M_+$ of all
positive elements $M.$ The following identity is crucial
 for the proof of the theorem:
$$
\Phi(\Delta(x)y)=-\Phi(x \Delta(y)),
$$
where $x, y$ are finite range operators from $LS(M).$

\begin{cor}\label{LSM}
Let $M$ be an arbitrary von Neumann algebra of type I$_\infty.$
 Then every $2$-local
derivation $\Delta: LS(M)\rightarrow LS(M)$
 is a
derivation.
\end{cor}

\section*{Acknowledgements}

The authors are indebted to the referee for valuable suggestions
and remarks.

\end{document}